\documentclass[twocolumn,letterpaper, 10 pt,conference]{ieeeconf}

\IEEEoverridecommandlockouts

\usepackage{amssymb}
\usepackage{amsmath}
\usepackage{color}
\usepackage{graphics}
\usepackage{graphicx}
\usepackage{multirow}
\usepackage{cite}
\usepackage{algorithm,algorithmic}
\usepackage{savesym}
\savesymbol{AND}
\usepackage{epstopdf}
\usepackage{subcaption}
\usepackage{float}
\usepackage{tikz}
\usetikzlibrary{angles,quotes}
\usetikzlibrary{calc}
\usetikzlibrary{arrows,decorations.markings}
\usepackage{mathtools}

\newtheorem{theorem}{Theorem}
\newtheorem{remark}{Remark}
\newtheorem{lemma}{Lemma}
\newtheorem{corollary}{Corollary}
\newtheorem{definition}{Definition}
\newtheorem{proposition}{Proposition}

\newtheorem{assumption}{Assumption}

\title{\LARGE \bf Probabilistic guarantees on the objective value for the scenario approach via sensitivity analysis}

\author{Zheming Wang and Rapha\"el M. Jungers
	\thanks{The authors are with the ICTEAM Institute, UCLouvain, Louvain-la-Neuve,1348, Belgium.
		R Jungers is a  FNRS honorary Research Associate.
		This project has received funding from the European Research Council (ERC) under the European Union's Horizon 2020 research and innovation programme under grant agreement No 864017 - L2C.
		He is also supported by the Walloon Region and the Innoviris Foundation.
		Email addresses:  \{zheming.wang, raphael.jungers\}@uclouvain.be}}

\begin{document}

\maketitle

\begin{abstract}
	This paper is concerned with objective value performance of the scenario approach for robust convex optimization. A novel method is proposed to derive probabilistic bounds for the objective value from scenario programs with a finite number of samples. This method relies on a max-min reformulation and the concept of complexity of robust optimization problems. With additional continuity and regularity conditions, via sensitivity analysis, we also provide explicit bounds which outperform an existing result in the literature. To illustrate the improvements of our results, we also provide a numerical example.
\end{abstract}

\section{Introduction}
In real-world applications, optimization problems with uncertainties often arise due to the lack of knowledge on system parameters and measurement/estimation errors. Classical approaches to tackle uncertain optimization problems include stochastic optimization in which uncertainties are treated as random variables \cite{BOO:P95,ART:RS03,BOO:BL11} and robust optimization which seeks to obtain robust solutions against all realizations of uncertainties \cite{ART:BN98,ART:BN08,BOO:BEN09,ART:BBC11}. An alternative approach for handling uncertainties is called the scenario approach or scenario optimization, which is somehow at the intersection of stochastic optimization and robust optimization as it is able to provide solutions to robust optimization and chance-constrained stochastic optimization problems in a probabilistic sense, see, e.g., some early works \cite{ART:CC05,ART:CC06}. We refer to \cite{BOO:CG18,ART:CCG21} for a comprehensive view on this subject.

The scenario approach is typically focused on constrained optimization problems with robustness constraints or chance constraints. By randomly sampling a finite set of constraints, a scenario program is formulated and a sample-based solution can be obtained. Chance-constrained theorems are then derived for this sample-based solution in terms of the measure of the violating subset of the uncertain constraints for both convex \cite{ART:CC05,ART:CC06,ART:CG08,ART:C10,ART:CG18,ART:GC19} and non-convex  \cite{ART:CGR18,ART:GC21} cases by using the concepts of support/essential constraints. In particular, the bounds in \cite{ART:CG08,ART:C10} have been proved to be tight for a special class of uncertain convex programs called fully supported problems.  While these theoretic results bring a lot of new perspectives into uncertain optimization problems, it takes further developments to apply the scenario approach to practical applications. For instance, many real-life physical properties of interest are directly related with the objective value of the underlying optimization problem, which means that the aforementioned chance-constrained theorems alone do not allow us to learn such properties. For this reason, formal analysis of the scenario approach is needed with respect to objective value performance.  To the best of our knowledge, the first work providing such results is \cite{ART:ESL14} under some continuity and regularity conditions on the basis of the chance-constrained theorems in \cite{ART:CG08,ART:C10}. Though the results in \cite{ART:ESL14} are built on an intermediate result that is tight in terms of the measure of the violating subset of the uncertain constraints (as proved in \cite{ART:CG08,ART:C10}), they are not necessarily tight in terms of the objective value. In this paper, we show that they are indeed not tight by providing improved bounds under the same conditions. In addition, the chance-constrained theorems in \cite{ART:CG08,ART:C10} require a non-degeneracy assumption. We show that such an assumption can be released when one is only concerned with the objective value.

This paper is motivated by recent progress on data-driven analysis of complicated systems, see, e.g., \cite{INP:KDSA14,INP:KQHT16,ART:BTSS20,ART:MXM20,ART:KBJT19,ART:WJ20,INP:BJW21,INP:RWJ21,ART:WJ21b}. In particular, for switched linear systems \cite{ART:KBJT19,ART:WJ20,INP:BJW21,INP:RWJ21,ART:WJ21b}, the scenario approach has been applied to the computation of the growth rate of the system using different types of Lyapunov functions after reformulating the stability analysis problem as a robust optimization problem in which the Lyapunov inequality is satisfied for all realizations of the trajectory. Different techniques are then proposed in these works to derive probabilistic guarantees on the computed growth rate. It appears that some of the results can be generalized from the stability analysis problem to general uncertain optimization problems. In \cite{INP:RWJ21}, it has been found that the solution to the stability analysis problem with an infinite number of constraints is actually defined by a finite set of support constraints whose cardinality is bounded by the dimension of the decision variable. This generalizes the definition of support constraints from scenario programs to robust optimization problems. This observation is then leveraged to derive probabilistic stability certificates using a sensitivity analysis argument on these support constraints.  In this paper, we extend this sensitivity analysis method to general robust optimization. We formally call the minimal number of such support constraints the \emph{complexity} of the robust optimization problem.

The contributions of this paper are two-fold: First, we present a general method of deriving probabilistic guarantees on the objective value performance of the scenario approach for robust convex programs. Second, we show theoretic improvements of our results as compared to \cite{ART:ESL14} under the same regularity conditions.

The rest of the paper is organized as follows. This section ends with the notation, followed by Section \ref{sec:pre} on the review of preliminary results on the scenario approach. Section \ref{sec:objective} presents a general method of deriving objective value performance guarantees for the scenario approach using a max-min reformulation. In Section \ref{sec:sensitivity}, under additional regularity conditions, we show explicit objective value performance bounds via sensitivity analysis.  Finally, an illustrative example is given in Section \ref{sec:example}.

\textbf{Notation}.
The non-negative real number set and the non-negative integer set are denoted by $\mathbb{R}^+$ and   $\mathbb{Z}^+$ respectively. Given any $r\in \mathbb{R}^+$ and $x\in \mathbb{R}^n$, $\mathbb{B}_r(x)$ denotes the open ball in $\mathbb{R}^n$ centered at $x$ with radius $r$. For a square matrix $P$, $P\succ(\succeq)~0$ means that $P$ is symmetric and positive definite (semi-definite), $\textrm{tr}(P)$ denotes the trace of $P$ and $\det(P)$ denotes the determinant of $P$. For a symmetric matrix $P$, we denote by $\lambda_{\max}(P)$ and $\lambda_{\min}(P)$ the largest and smallest eigenvalues of $P$ respectively. For any matrix $P\succ 0$, let $\kappa(P)\coloneqq \lambda_{\max}(P)/\lambda_{\min}(P)$ and $\chi(P)\coloneqq \sqrt{\det(P)/\lambda_{\min}(P)^n}$. For any $p\ge 1$, the $\ell_p$ norm of a vector $x\in \mathbb{R}^n$ is $\|x\|_p$ ($\|x\|$ is the $\ell_2$ norm by default). For a set $S$, $\textrm{cl}(S)$ denote the closure of $S$.

\section{Preliminaries}\label{sec:pre}
We consider the following robust convex optimization problem (RCP):
\begin{subequations}\label{eqn:cxfxd}
\begin{align}
 J^\star:=&\inf_{x\in X} c^\top x\\
\textrm{s.t. }& f(x,\delta) \le 0, \forall \delta \in \Delta
\end{align}
\end{subequations}
where $x\in  \mathbb{R}^d $ is the decision variable, $\delta$ is the uncertain parameter contained in the uncertain set $\Delta \subseteq \mathbb{R}^m$ with the uncertain constraint $f(x,\delta) \le 0$, and $X\subseteq \mathbb{R}^d$ is some static constraint set. We would like to mention that, while there is only one uncertain constraint in (\ref{eqn:cxfxd}), the results and discussions in the paper extend trivially to cases with multiple uncertain constraints. The following assumptions are made.

\begin{assumption}\label{ass:continuous}
	$f: X \times \Delta \rightarrow \mathbb{R}$ is lower semi-continuous (l.s.c.) in $X \times \Delta$ and $X\subseteq \mathbb{R}^d$ is compact.
\end{assumption}

\begin{assumption}\label{ass:feasible}
	There exists $x\in X$ such that $f(x,\delta)\le 0$ for any $\delta \in \Delta$.
\end{assumption}

\begin{assumption}\label{ass:convex}
	$X\subseteq \mathbb{R}^d$ is convex and $f: X \times \Delta \rightarrow \mathbb{R}$ is convex in $x\in X$ for any $\delta \in \Delta$.
\end{assumption}


%



Let $(\Delta,\mathcal{B}(\Delta),\mathbb{P})$ be a probability space where $\Delta$ is a metric space, $\mathcal{B}(\Delta)$ is the Borel $\sigma$-algebra of $\Delta$ and $\mathbb{P}: \mathcal{B}(\Delta) \rightarrow [0,1]$ is the probability measure.

In the framework of the scenario approach (also called scenario optimization) \cite{ART:CC05,ART:CC06,ART:CCG21}, the solution of the robust optimization problem in (\ref{eqn:cxfxd}) is approximated with a finite number of samples. More precisely, we randomly generate $N\in \mathbb{Z}^+$ samples, denoted by $(\delta_1,\delta_2,\cdots, \delta_N) \in \Delta^N$, from some probability measure with the support $\Delta$, and formulate the following scenario convex program (SCP)
\begin{subequations}\label{eqn:deltaD}
\begin{align}
\mathcal{P}(\omega_N): \quad J(\bar{\omega}_N)\coloneqq&\min_{x\in X} c^\top x\\
\textrm{s.t. }& f(x,\delta) \le 0, \forall \delta \in \omega_N.
\end{align}
\end{subequations}
where $\omega_N : = (\delta_1,\delta_2,\cdots, \delta_N)$, and $\bar{\omega}_N$ denote the unordered sample set, i.e.,  $\bar{\omega}_N \coloneqq \{\delta_1,\delta_2,\cdots, \delta_N\}$, and $J: 2^{\Delta} \rightarrow \mathbb{R}$ with $J(\Delta) = J^{\star}$. 

Chance-constrained guarantees on the solution of the sampled problem are derived in \cite{ART:CG08,ART:C10} for convex problems in terms of the violation probability, depending on the dimension of the decision variable and the number of samples, as shown below.

\begin{theorem}[{\cite[Theorem 1]{ART:CG08},\cite[Theorem 3.3]{ART:C10}}] \label{thm:chance}
For any $N\in \mathbb{Z}^+$ with $N>d$, suppose $\omega$ are $N$ independent and identically distributed (i.i.d.) samples drawn according to the probability measure $\mathbb{P}$. Consider $\mathcal{P}(\omega)$ as defined in (\ref{eqn:deltaD}), let us assume that $\mathcal{P}(\omega)$ admits a unique optimal solution. Then, with  Assumptions \ref{ass:continuous}--\ref{ass:convex} and additional conditions (see Remark \ref{rem:con}), for any $\epsilon \in (0,1)$,
\begin{align}\label{eqn:ccbound}
\mathbb{P}^{N} \{\omega \in \Delta^N: \mathcal{V}(\bar{\omega})> \epsilon \} 
\le \Phi_c(\epsilon;d,N)
\end{align}
where $\mathcal{V}(\bar{\omega})\coloneqq\mathbb{P}\{\delta: J(\bar{\omega}\cup \{\delta\})>J(\omega)\}$, $\bar{\omega}$ is the set of the elements in $\omega$, $J(\bar{\omega})$ is defined as in (\ref{eqn:deltaD}), $\mathbb{P}^N$ is the probability measure in the $N$-Cartesian product of $\Delta$, and
\begin{align}\label{eqn:phic}
	\Phi_c(\epsilon;k,N) &\coloneqq \sum_{i=0}^{k-1} {N \choose i} \epsilon^i(1-\epsilon)^{N-i}, \forall k\ge 1.
\end{align}
\end{theorem}

In fact, the bound in (\ref{eqn:ccbound}) is proved to be tight for a special class of convex programs called fully supported programs, see \cite{ART:CG08,ART:C10} for details. For non-convex optimization problems, similar results are derived in \cite{ART:CGR18,ART:GC21}. While these probabilistic guarantees for both the convex and non-convex cases provide a fundamental understanding on the solution of the scenario program \eqref{eqn:deltaD}, they alone do not allow to bound the objective value of the original robust optimization.

\begin{remark}\label{rem:con}
The bound in (\ref{eqn:ccbound}) is derived by \cite{ART:CG08} and \cite{ART:C10} based on different conditions. More precisely, \cite{ART:CG08} requires that the feasible domain of $\mathcal{P}(\omega)$ has a nonempty interior for any $\omega$; \cite{ART:C10} assumes that $\mathcal{P}(\omega)$ is nondegenerate with probability one, see Definition 2.7 in \cite{ART:C10} for the definition of degenerate problems.
\end{remark}

In \cite{ART:ESL14}, theoretical links between $\mathcal{P}(\omega)$ and the original robust optimization problem RCP in (\ref{eqn:cxfxd}) are further established in terms of the objective value in a statistical sense under the assumptions in Theorem \ref{thm:chance} and Slater's condition, i.e., there exists $x\in X$ such that $\sup_{\delta\in \Delta} f(x,\delta) <0$. We briefly summarize the idea in \cite{ART:ESL14} as follows. First, it is shown that the optimizer of SCP is a feasible solution to the following chance-constrained program (CCP) with some confidence level
\begin{subequations}\label{eqn:ccpep}
	\begin{align}
	&\min_{x\in X} c^\top x\\
		\textrm{s.t. } & \mathbb{P}\{f(x,\delta) \le 0 \} \ge 1-\epsilon
	\end{align}
\end{subequations}
where $\epsilon\in (0,1)$.
It is then shown that the feasible set of the program CCP is contained in the feasible set of a relaxed RCP in the form of 
\begin{subequations}\label{eqn:cxfxdrelax}
	\begin{align}
		&\inf_{x\in X} c^\top x\\
		\textrm{s.t. }& f(x,\delta) \le \xi, \forall \delta \in \Delta
	\end{align}
\end{subequations}
where $\xi \ge 0$. Finally, with the relations above, probabilistic objective bounds for the scenario approach can be derived from the relaxed RCP under Lipschitz continuity and other regularity assumptions. We call this procedure a two-step approach in the sense that it relies on the link from SCP to CCP (via the chance-constrained Theorem \ref{thm:chance}) to establish probabilistic objective value bounds. As we will show in the sequel, the step from SCP to CCP brings conservatism into the resulting bound.

In this paper, we revisit the problem of deriving probabilistic objective value bounds for the scenario approach from a sensitivity analysis perspective in the context of robust optimization. We develop a method that skips the step from SCP to CCP and enables us to directly establish probabilistic objective value bounds for SCP using a perturbation argument. Thus, our approach does not require the conditions in Remark \ref{rem:con}.  We believe that this latter fact constitutes an important improvement because these conditions can be violated in some control problems, see \cite{INP:BJW21} for such an example.

\section{Probabilistic objective value performance}\label{sec:objective}
This section presents a general method deriving probabilistic objective value bounds of the scenario approach for robust optimization problems. This method adopts a max-min reformulation of the original robust optimization problem (\ref{eqn:cxfxd}).

\subsection{Continuity properties}
We first show some preliminary continuity results on parametric optimization, see, e.g., \cite{BOO:BGKKT83,BOO:RW98,ART:S18} for a comprehensive view. Given any $N\in \mathbb{Z}^+$, for any $\omega_N\coloneqq(\delta_1,\delta_2,\cdots, \delta_N)\in \Delta^N$, for ease of discussion, let us define
\begin{align}
	g_N(\omega_N) &:= J(\{\delta_1,\delta_2,\cdots, \delta_N\}), \label{eqn:gddelta}\\
	F_N(\omega_N) &:=\{x\in X: f(x,\delta_i) \le 0, i=1,2,\cdots, N\}, \label{eqn:Fddelta}\\
	F^o_N(\omega_N) &:=\{x\in X: f(x,\delta_i) < 0, i=1,2,\cdots, N\},  \label{eqn:Foddelta}\\
	\phi_N(\omega_N) &:= \{x\in F_N(\omega_N): c^\top x =  g_N(\omega_N)\}, \label{eqn:phiddelta}
\end{align}
Under Assumption \ref{ass:continuous}, a semi-continuity property is stated below. 
\begin{lemma}\label{lem:lsc}
	Consider the RCP in \eqref{eqn:cxfxd} and the SCP in \eqref{eqn:deltaD}. Given any $N\in \mathbb{Z}^+$, for any $\omega\in \Delta^N$, let $g_N(\omega)$ be defined as in (\ref{eqn:gddelta}). 
	Suppose Assumptions \ref{ass:continuous}\& \ref{ass:feasible}  hold, then $g_N(\omega)$ is l.s.c. in $\Delta^N$.
\end{lemma}
Proof: We prove the lemma by using standard results on parametric optimization. With Assumptions \ref{ass:continuous}\& \ref{ass:feasible},  from Theorem 4.2.1 in \cite{BOO:BGKKT83}, to show the l.s.c. property of $g_N(\omega)$, we only need to show that the point-to-set mapping $F_N: \Delta^N \rightarrow 2^{X}$, as defined in (\ref{eqn:Fddelta}), is closed at any $\omega \coloneqq (\delta_1,\delta_2,\cdots, \delta_N)\in \Delta^N$ in the sense that, for each pair of sequences $\{\omega^t\coloneqq (\delta_1^t,\delta_2^t,\cdots, \delta_N^t), t\in \mathbb{Z}^+\} \subset \Delta^N$ and $\{x^t,t\in \mathbb{Z}^+\} \subset X$ with the properties $\omega^t \rightarrow \omega$, $x^t \in F_N(\omega^t)$, and $x^t \rightarrow x$, it follows that $x\in F_N(\omega)$. From the fact that $f: X \times \Delta \rightarrow \mathbb{R}$ is l.s.c. in $X \times \Delta$, $f(x,\delta_i) \le \lim\inf_{(x^t,\delta^t_i)\rightarrow (x,\delta_i)}f(x^t,\delta_i^t) \le 0$ for any $i=1,2,\cdots, N$, where the second inequality is due to the property that $x^t \in F_N(\omega^t)$. Hence, we conclude $x\in F_N(\omega)$. A similar proof can also be derived from Lemma 5.3 in \cite{ART:S18}. $\Box$

\subsection{Max-min reformulation}


With the continuity property above, we then define, for any $N\in \mathbb{Z}^+$, 
\begin{align}\label{eqn:gstd}
g^*_N := \sup_{\omega \in \Delta^N} g_N(\omega) = \sup_{\omega \in \Delta^N} \inf_{x\in X} c^\top x+ \chi_{S_N} (x,\omega),
\end{align}
where $S_N:= \{(x,\delta_1,\delta_2,\cdots, \delta_N): f(x,\delta_i)\le 0, i=1,2\cdots, N)\}$ and $\chi_{S_N}: X \times \Delta^N \rightarrow \{0,+\infty\}$ is the characteristic function of $S_N$ defined as follows: $\chi_{S_N}(x,\omega) = 0$ if $(x,\omega)\in S_N$ and  $ = +\infty$ everywhere else.

The next lemma states the equivalence between the RCP in (\ref{eqn:cxfxd}) and the max-min formulation as defined in (\ref{eqn:gstd}).

\begin{lemma}\label{lem:JG}
Consider Problems (\ref{eqn:cxfxd}) and (\ref{eqn:gstd}) for any $N\in \mathbb{Z}^+$ with the optimal values $J^\star$ and $g^*_N$ respectively. Suppose Assumptions \ref{ass:continuous} -- \ref{ass:convex} hold, then it holds that $J^\star=g^*_N$ for any $N\ge d$; 
\end{lemma}
Proof: From the definitions in (\ref{eqn:cxfxd}) and (\ref{eqn:gstd}), it holds that $g_N(\omega)\le J^\star$ for any $\omega \in \Delta^N$, which implies that $g^*_N\le J^\star$ for any $N\in \mathbb{Z}^+$. As $\{g^*_N\}_{N\in \mathbb{Z}^+}$ is a monotonically increasing sequence, bounded from above by $J^\star$, we only need to show that $g^*_d =J^\star$. First, we show that $g^*_{d+1}=J^\star$. The optimality of $g^*_{d+1}$ in (\ref{eqn:gstd}) implies that the set $ \mathcal{X}_{(\delta_1,\delta_2,\cdots, \delta_{d+1})} \coloneqq \{x\in X: c^\top x \le g^*_d, f(x,\delta_i) \le 0, i=1,2,\cdots, d+1 \}$ is nonempty for any $(\delta_1,\delta_2,\cdots, \delta_{d+1})\in \Delta^{d+1}$. Since $f(x,\delta)$ is l.s.c. and $X$ is a compact convex set, $\mathcal{X}_{(\delta_1,\delta_2,\cdots, \delta_{d+1})}$ is also a compact convex set. Hence, from Helly's theorem (see, e.g., Theorem 21.6 in \cite{BOO:R70}), $\{x\in X: c^\top x \le g^*_{d+1}, f(x,\delta) \le 0, \delta\in \Delta\}$ is nonempty, which means that $J^\star\le g^*_{d+1}$ and thus $g^*_{d+1}=J^\star$. Now, we show that $g^*_d = g^*_{d+1}$. For any  $\omega_{d+1} \coloneqq(\delta_1,\delta_2,\cdots, \delta_{d+1})\in \Delta^{d+1}$, from  Theorem 2 in \cite{ART:CC05} (or Theorem 3 in \cite{ART:CC06}), there exists $i$ such that $J(\bar{\omega}_{d+1}) = J(\bar{\omega}_{d+1}\setminus \{\delta_i\})$, where $\bar{\omega}_{d+1}\coloneqq \{\delta_1,\delta_2,\cdots, \delta_{d+1}\}$. As a result, $g^*_d = g^*_{d+1}$, which implies that $g^*_d = J^\star$. $\Box$

In view of the result in Lemma \ref{lem:JG}, we define the concept of \emph{complexity} for robust optimization problems with infinitely many constraints, following the definition of \emph{complexity} for scenario programs with a finite number of constraints in \cite{ART:GC19}.

\begin{definition}\label{def:helly}
	Consider the RCP in the form of (\ref{eqn:cxfxd}), the \emph{complexity} of RCP is the smallest integer $k\in \mathbb{Z}^+$ such that $g^*_k=J^\star$.
\end{definition}

Similar to the \emph{complexity} of SCP in \cite{ART:GC19}, the \emph{complexity} of RCP is intuitively the minimal number of samples needed to obtain the true objective value $J^\star$ with random sampling. With this formal definition, Lemma \ref{lem:JG} basically means that the \emph{complexity}  of RCP is bounded by the number of decision variables under Assumptions \ref{ass:continuous} -- \ref{ass:convex}.

\begin{remark}
When RCP is not feasible, i.e., Assumption \ref{ass:feasible} does not hold or $J^\star=+\infty$, from Helly's theorem, it still holds that $J^\star=g^*_{d+1}$ but not $J^\star=g^*_d$ in Lemma \ref{lem:JG}. Thus, the \emph{complexity} of RCP is bounded by $d+1$ in any case.
\end{remark}

\subsection{Uniform level-set bounds}\label{sec:ULB}
To establish probabilistic objective value performance guarantees, inspired from \cite{ART:KT12,ART:ESL14}, we define the following tail probability for any $N\in \mathbb{Z}^+$
\begin{align}\label{eqn:palpha}
 p_N(\alpha):= \mathbb{P}^N\{\omega\in \Delta^N: g^*_N - \alpha \le g_N(\omega)\}, \forall \alpha \ge 0.
\end{align}
From the semi-continuity property of $g_N$ in Lemma \ref{lem:lsc}, $ p_N(\alpha)$ is well defined and monotonically increasing. In particular, when $N\ge d$, from Lemma \ref{lem:JG}, $ p_N$ defined as in (\ref{eqn:palpha}) becomes
\begin{align}
 p_N(\alpha)= \mathbb{P}^N \{\omega\in \Delta^N: J^\star - \alpha \le g_N(\omega)\}, \forall \alpha \ge 0.
\end{align}
The continuity property of this tail probability function is stated in the following lemma.

\begin{lemma}\label{lem:semiconupper}
Consider the RCP \eqref{eqn:cxfxd} and the SCP \eqref{eqn:deltaD}. Suppose Assumptions \ref{ass:continuous} - \ref{ass:convex} hold. For any $N\in \mathbb{Z}^+$, let the tail probability function $ p_N: \mathbb{R}^+ \rightarrow [0,1]$ be defined as in (\ref{eqn:palpha}). Then, the following property holds:  i) There exists a finite $C\in \mathbb{R}^+$ such that $p_D(\alpha)=1$ for any $\alpha \ge C$. ii) $ p_N$ is upper semi-continuous (u.s.c.) at any $\alpha \in \mathbb{R}^+$.
\end{lemma}
Proof: i) Let $\tilde{J}:= \min_{x\in X} c^\top x$. As a result, $g_N(\omega) \ge \tilde{J}$ for any $z\in \Delta$. Hence, for any $\alpha \ge p_N^*-\tilde{J}$, $p_N(\alpha) = 1$.\\
ii) It is easy to verify that $ p_N(\alpha)$ is monotonically nondecreasing. Now, we show that for any $\bar{\alpha} \in \mathbb{R}^+$, it holds that $p_N(\bar{\alpha}) \ge \limsup_{\alpha \rightarrow \bar{\alpha}} p_N(\alpha)$, i.e., $\lim_{i\rightarrow \infty} \alpha_i = \bar{\alpha}$ implies $p_N(\bar{\alpha}) \ge \limsup_{i \rightarrow \infty} p_N(\alpha_i)$. From the monotonicity of $p_N(\alpha)$, we only need to consider the case that $\alpha_{i+1} \le \alpha_{i}$ for all $i\in \mathbb{Z}^+$. We then define: $S_i=\{\omega\in \Delta^N: J^\star - \alpha_i \le g_N(\omega))\}$ for  $i\in \mathbb{Z}^+$. Thus, $\limsup_{i \rightarrow \infty} p_N(\alpha_i) = \lim_{i \rightarrow \infty} \mathbb{P}^N\{S_i\} =  \mathbb{P}^N\{\lim_{i \rightarrow \infty} S_i\} = p_N(\bar{\alpha})$ since $S_i\supseteq S_{i+1} $ and $\mathbb{P}^N\{S_i\}$ is bounded for any $i\in \mathbb{Z}^+$. This proves upper semi-continuity of $ p_N$. $\Box$


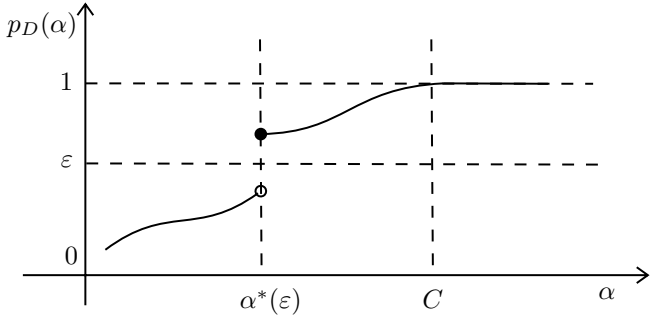
\begin{figure}

\tikzset{every picture/.style={line width=0.75pt}} 
\centering
\begin{tikzpicture}[x=0.6pt,y=0.6pt,yscale=-1,xscale=1]

\draw  (120,281) -- (514,281)(159.4,110) -- (159.4,300) (507,276) -- (514,281) -- (507,286) (154.4,117) -- (159.4,110) -- (164.4,117)  ;
\draw    (172,265) .. controls (211.4,235.45) and (229.46,257.32) .. (268.21,229.31) ;
\draw [shift={(270,228)}, rotate = 323.13] [color={rgb, 255:red, 0; green, 0; blue, 0 }  ][line width=0.75]      (0, 0) circle [x radius= 3.35, y radius= 3.35]   ;
\draw    (270,192) .. controls (324,191) and (326,163) .. (385,160) ;
\draw [shift={(270,192)}, rotate = 358.94] [color={rgb, 255:red, 0; green, 0; blue, 0 }  ][fill={rgb, 255:red, 0; green, 0; blue, 0 }  ][line width=0.75]      (0, 0) circle [x radius= 3.35, y radius= 3.35]   ;
\draw    (385,160) -- (451,160.25) ;
\draw  [dash pattern={on 4.5pt off 4.5pt}]  (159.5,160) -- (479.5,160.25) ;
\draw  [dash pattern={on 4.5pt off 4.5pt}]  (160,210.5) -- (483,211.25) ;
\draw  [dash pattern={on 4.5pt off 4.5pt}]  (269.5,132.25) -- (270.5,280.75) ;
\draw  [dash pattern={on 4.5pt off 4.5pt}]  (377.5,132.25) -- (378.5,280.75) ;

\draw (108.5,112.9) node [anchor=north west][inner sep=0.75pt]    {$p_{D}( \alpha )$};
\draw (141.5,151.4) node [anchor=north west][inner sep=0.75pt]    {$1$};
\draw (481,286.9) node [anchor=north west][inner sep=0.75pt]    {$\alpha $};
\draw (144.5,261.4) node [anchor=north west][inner sep=0.75pt]    {$0$};
\draw (142,203.4) node [anchor=north west][inner sep=0.75pt]    {$\varepsilon $};
\draw (255,287.4) node [anchor=north west][inner sep=0.75pt]    {$\alpha ^{*}( \varepsilon )$};
\draw (370,288.4) node [anchor=north west][inner sep=0.75pt]    {$C$};

\end{tikzpicture}

\caption{Illustration of the tail probability function and the optimal ULB: Given $\varepsilon\in (0,1)$, $\alpha ^{*}( \varepsilon )$ is the supremum of $\alpha$ satisfying $ p_N(\alpha)\le \varepsilon$. }\label{fig:tail}
\end{figure}

We then call $\tilde{\alpha} : (0,1) \rightarrow \mathbb{R}$ a uniform level-set bound (ULB) if it satisfies
\begin{align}\label{eqn:pDepsilon}
p_N(\tilde{\alpha}(\varepsilon) ) \ge \varepsilon, \forall \varepsilon\in (0,1).
\end{align}
The optimal ULB is defined as follows:
\begin{align}\label{eqn:alphastar}
\alpha^*(\varepsilon)\coloneqq \sup_{\alpha\ge 0} \alpha: p_N(\alpha)\le \varepsilon, \forall \varepsilon\in (0,1).
\end{align}
With the discussions above, we arrive at the main result of this section.

\begin{proposition}\label{prop:boundd}
Consider the RCP in \eqref{eqn:cxfxd} and the SCP in \eqref{eqn:deltaD}. Suppose Assumptions \ref{ass:continuous} -- \ref{ass:convex} hold. Given $N (N\ge d)$ i.i.d. samples, denoted by $\omega_N$,  drawn according to the probability measure $\mathbb{P}$, consider Problem (\ref{eqn:deltaD}) with the optimal value $g_N(\omega_N)$ as defined in \eqref{eqn:gddelta}. For any $\varepsilon \in (0,1)$, let $\alpha^*(\varepsilon)$ be defined as in  (\ref{eqn:alphastar}). Then, it holds that
\begin{align}\label{eqn:JstarJalphaep}
\mathbb{P}^{N}\{\omega_N \in \Delta^{N}: J^\star \le g_N(\omega_N)+ \alpha^*(\varepsilon)  \} \ge \varepsilon
\end{align}
where $J^\star$ is given in (\ref{eqn:cxfxd}).
\end{proposition}
Proof: With Lemma \ref{lem:JG} at hand, we only need to show $p_N(\alpha^*(\epsilon)) \ge \varepsilon$ for any $\varepsilon \in (0,1)$ as $J^\star = g_N^*$. We consider two cases: i) the point at the optimal is continuous; ii)  the  point at the optimal is not continuous. In the first case, from the monotonicity of $p_N(\alpha)$, $p_N(\alpha)\le \varepsilon$ is active, meaning that $p_N(\alpha^*(\varepsilon))=\varepsilon$. In the second case, from the upper semi-continuity shown in Lemma \ref{lem:semiconupper}, $p_N(\alpha^*(\varepsilon)) \ge \varepsilon$. Thus, we conclude that  $p_N(\alpha^*(\varepsilon)) \ge \varepsilon$. $\Box$

While this is an elementary result, following from the bound on the \emph{complexity} of RCP and the u.s.c. of $p_N$, as illustrated in Figure \ref{fig:tail}, it provides a general way of bounding the original objective value $J^\star$. In practice, as computing the optimal ULB is not always possible, we turn to upper bounds of $\alpha^*(\varepsilon)$.


\begin{remark} From the proof of Lemma \ref{lem:semiconupper}, 
	when there is no measure concentration on the level sets of $g_N$, i.e., $\mathbb{P}^N(\omega\in \Delta^N: g_N(\omega) = u)=0$ for any $u\in \mathbb{R}$, it can be shown that $p_N(\cdot)$ is continuous. Thus, the optimal ULB can be defined by $p_N(\alpha) = \varepsilon$.
\end{remark}

The definition of the tail probability function in (\ref{eqn:palpha}) can be also extended to a functional space as follows: $\forall \alpha \in \mathcal{L}^+_N$,
\begin{align}
 p_N(\alpha)= \mathbb{P}^N \{\omega\in \Delta^N: J^\star - \alpha(\omega) \le g_N(\omega)\}.
\end{align}
where $\alpha: \Delta^N \rightarrow \mathbb{R}^+$ and $\mathcal{L}^+_N$ is the space of measurable positive real functions in $\Delta^N$. Similarly, probabilistic upper bounds on $J^\star$ can be obtained by finding $\alpha \in  \mathcal{L}^+_N$ such that $ p_N(\alpha) \ge \varepsilon$ for some given $\varepsilon \in (0,1)$.

\section{Explicit bounds via sensitivity analysis}\label{sec:sensitivity}
We now explicitly characterize ULBs as defined in (\ref{eqn:pDepsilon}) via sensitivity analysis. 

\subsection{Sensitivity analysis}
With an assumption on the continuity of $g_d: \Delta^d \rightarrow \mathbb{R}$, a deterministic relation between the robust optimization problem RCP and the scenario program SCP is established in the following lemma. We will provide sufficient conditions to ensure this assumption later in this section.
\begin{lemma}\label{lem:LgJ}
Suppose Assumptions \ref{ass:continuous} - \ref{ass:convex} hold. 
Let $g_d(\omega)$ be defined as in (\ref{eqn:gddelta}) for any $\omega\in \textrm{cl}(\Delta)^d$. Assume that there exists a constant $L_g$ such that, for any $\omega=(\delta_1,\delta_2,\cdots,\delta_d) \in \textrm{cl}(\Delta)^d$ and $\omega'=(\delta'_1,\delta'_2,\cdots,\delta'_d) \in \textrm{cl}(\Delta)^d$,  
\begin{align}\label{eqn:gdLg}
	\|g_d(\omega)-g_d(\omega')\| \le L_g \max_{i\in \{1,2,\cdots, d\}} \|\delta_i-\delta'_i\|.
\end{align}
Then,  the following results holds: (i)  There exists a set $\bar{\omega}_d^*\coloneqq \{\delta^*_1,\delta^*_2,\cdots, \delta^*_d\} \subset \textrm{cl}(\Delta)$ such that $J(\bar{\omega}_d^*) = J^\star$. (ii) Given any finite subset $\bar{\omega} \subseteq \textrm{cl}(\Delta)$, $J(\bar{\omega}) \ge J^\star - L_g \|\bar{\omega}_d^*\|_{\bar{\omega}} $.
\end{lemma}
Proof: (i) From Lemma \ref{lem:JG}, we know that $J^\star = g^*_d := \sup_{\omega \in \Delta^d} g_d(\omega)$. The condition in \eqref{eqn:gdLg} implies Lipschitz continuity of $g_d(\omega)$ in $\textrm{cl}(\Delta)^d$, which means that, according to Weierstrass's theorem, there exist $\omega\in \textrm{cl}(\Delta)^d$ such that $g^*_d = g(\omega)$.\\
(ii) Let us define the set $\bar{\omega}_d \coloneqq \{\delta_1,\delta_2,\cdots, \delta_d\} \subseteq \bar{\omega}$ as follows: $\delta_i =\arg\min_{\delta\in \bar{\omega}} \|\delta_i^*-\delta\|$ for any $i=1,2,\cdots, d$.
 Thus, from \eqref{eqn:gdLg},
\begin{align*}
J(\bar{\omega}_d^*)-J(\bar{\omega}) &\le  J(\bar{\omega}_d^*)-J(\bar{\omega}_d) = g_d(\omega_d^*)-g_d(\omega_d)\\
& \le L_g \max_{i\in \{1,2,\cdots, d\}} \|\delta^*_i-\delta_i\| =  L_g \|\bar{\omega}_d^*\|_{\bar{\omega}},
\end{align*}
where $\omega_d^*\coloneqq (\delta^*_1,\delta^*_2,\cdots, \delta^*_d)$ and $\omega_d \coloneqq (\delta_1,\delta_2,\cdots, \delta_d)$. $\Box$


When the uncertainty $\delta$ in RCP is a scalar, the condition \eqref{eqn:gdLg} is simply the Lipschitz continuity with respect to the $\ell_\infty$ norm. In order to obtain an explicit ULB, we also need a regularity condition on $\Delta$ with respect to the measure $\mathbb{P}$. We consider the same condition as in \cite{ART:ESL14}. 

\begin{assumption}\label{ass:regular}
For the probability space $(\Delta,\mathcal{B}(\Delta),\mathbb{P})$, there exists a strictly increasing function $\varphi: \mathbb{R}^+ \rightarrow [0,1]$ such that
\begin{align}
\forall \delta \in \textrm{cl}(\Delta), \forall r\in \mathbb{R}^+, \mathbb{P}\{\mathbb{B}_r(\delta)\cap \Delta\} \ge \varphi(r), 
\end{align}
where $\mathbb{B}_r(\delta)$ is an open ball centered at $\delta$ with radius $r$.
\end{assumption}

\begin{remark}
	Inspired by \cite{ART:CR04}, we can also consider a more explicit regularity condition. There exist $c>0$ and $\bar{r}>0$ such that
	\begin{align}
		\mathbb{P}\{\mathbb{B}_r(\delta)\cap \Delta\} \ge c\mu(\mathbb{B}_r(\delta)), \forall \delta \in \textrm{cl}(\Delta), 0 < r \le \bar{r}
	\end{align}
	where $\mu(\cdot)$ denotes the Lebesque measure.
\end{remark}

The following lemma shows an elementary probability result.

\begin{lemma}\label{lem:Phi}
Given the probability space $(\Delta,\mathcal{B}(\Delta),\mathbb{P})$, for any $k\in \mathbb{Z}^+$ and $\epsilon \in [0,1]$,  consider $k$ subsets $\{B_1,B_2,\cdots,B_k\}$ of $\Delta$ with $\mathbb{P}\{B_i\} = \epsilon$ for all $i=1,2,\cdots, k$.  For any $N\in \mathbb{Z}^+$,  suppose $\omega$ is a vector of $N$ independent and identically distributed (i.i.d.) samples drawn according to the probability measure $\mathbb{P}$. Then,
\begin{align}\label{eqn:BPhi}
\mathbb{P}^N\{\omega\in \Delta^N: B_i \cap \bar{\omega} \not=\emptyset, \forall i\} \ge  1-\Phi_a(\epsilon;k,N)
\end{align}
where $\bar{\omega}$ is the set of elements in $\omega$ and
\begin{align}\label{eqn:phia}
	 \Phi_a(\epsilon;k,N) &\coloneqq k(1-\epsilon)^N.
\end{align}
%
\end{lemma}
Proof: For any $i$, let us define $S_i \coloneqq \{\omega\in \Delta^N: B_i \cap \bar{\omega} = \emptyset \}$, i.e.,  the set of $\omega$ such that none of the points in $\omega$ fails into $B_i$. As  $\mathbb{P}\{B_i\} = \epsilon$, we know that $\mathbb{P}^N\{S_i\} = (1-\epsilon)^N$ for any $i$. Hence, $\mathbb{P}^N\{ \cup_{i}^k S_i\} \le k(1-\epsilon)^N$. Note that $\cup_{i}^k S_i$ is the set $\omega$ such that there exists at least $i$ such that $B_i \cap \bar{\omega} = \emptyset$. Hence, $\mathbb{P}^N\{\bar{\omega}\in \Delta^N: B_i \cap \bar{\omega} \not=\emptyset, \forall i\} = 1-\mathbb{P}^N\{ \cup_{i}^k S_i\} \ge 1-k(1-\epsilon)^N$. $\Box$

In fact, a tighter bound can be derived when $\epsilon$ in Lemma \ref{lem:Phi} is sufficiently small. However, the proof is too long to be included in the paper. We only comment on this in the following remark.

\begin{remark}\label{rem:phi}
	With the same conditions as Lemma \ref{lem:Phi}, when $\epsilon \le 1/k$, it can be verified that
	\begin{align}\label{eqn:PhiN}
		\mathbb{P}^N\{\omega\in \Delta^N: B_i \cap \bar{\omega} \not=\emptyset, \forall i\} \ge 1- \Phi(\epsilon;k,N)
	\end{align}
where
	\begin{align}\label{eqn:Phi}
		\Phi(\epsilon;k,N):=\sum_{i=1}^{k}(-1)^{i-1} \sum_{j=0}^{k-i} { {i-1+j}\choose{j}} \left(1-i\epsilon \right)^N.
	\end{align}
It can be also shown that $\Phi(\epsilon;k,N) \le k(1-\epsilon)^N$.
\end{remark}

\begin{figure}
	\centering
	\includegraphics[width=0.9\linewidth]{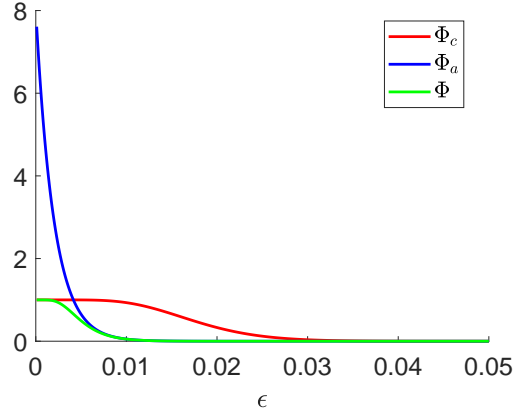}
	\caption{Comparison of different probability bounds with $k=8$ and $N=500$.}\label{fig:phi}
\end{figure}

We compare the bounds in (\ref{eqn:phic}) (\ref{eqn:phia}), and (\ref{eqn:BPhi})  in Figure \ref{fig:phi}. It can be seen that $\Phi(\epsilon;k,N)$ is the best among these three bounds. As $\epsilon$ increases, $\Phi_a(\epsilon;k,N)$ becomes close to $\Phi(\epsilon;k,N)$ and outperforms $\Phi_c(\epsilon;k,N)$ eventually.

Based the lemmas above, we then present the main result of this section.
\begin{theorem}\label{thm:gNJstar}
Consider Problem (\ref{eqn:cxfxd}) with the optimal objective value $J^\star$, suppose Assumptions \ref{ass:continuous} -- \ref{ass:regular} hold. For any $N\in \mathbb{Z}^+$ with $N \ge d$, 
let $g_N(\omega)$ be defined as in (\ref{eqn:gddelta}) for any $\omega\in \textrm{cl}(\Delta)^N$. Assume that there exists a constant $L_g$ such that \eqref{eqn:gdLg} is satisfied in $\textrm{cl}(\Delta)^d$. Let $\omega_N$ be a vector of $N$ i.i.d. samples drawn according to the probability measure $\mathbb{P}$. Then, for any $\epsilon\in (0,1)$, it holds that
\begin{align}\label{eqn:PNomegaJstargPhi}
&\mathbb{P}^N\{\omega_N \in \Delta^{N}: J^\star \le g_N(\omega_N)+L_{g} \varphi^{-1}(\epsilon)\} \nonumber\\
&\qquad \ge 1 - \Phi_a(\epsilon;d,N) = 1-d(1-\epsilon)^N.
\end{align}
\end{theorem}
Proof: From property (i) of Lemma \ref{lem:LgJ}, there exists a set $\bar{\omega}_d^*\coloneqq \{\delta^*_1,\delta^*_2,\cdots, \delta^*_d\} \subset \textrm{cl}(\Delta)$ such that $J(\bar{\omega}_d^*) = J^\star$. For any $\epsilon\in (0,1)$, let $r_i$ be such that $\mathbb{P}\{\mathbb{B}_{r_i}(\delta^*_i)\cap \Delta\} = \epsilon$ for any $i=1,2,\cdots, d$. Under Assumption \ref{ass:regular}, we can see that $r_i \le \varphi^{-1}(\epsilon)$ for any $i=1,2,\cdots, d$. For convenience, we define $B_i = \mathbb{B}_{r_i}(\delta^*_i)\cap \Delta$ for any $i=1,2,\cdots, d$. With Lemma \ref{lem:Phi}, we know that $\mathbb{P}^N\{\omega_N\in \Delta^N: B_i \cap \bar{\omega}_N \not=\emptyset, \forall i=1,2,\cdots, d\} \ge  1-d(1-\epsilon)^N$, where $\bar{\omega}_N$ is the set of the elements in the stacked vector $\omega_N$. Since $B_i \cap \bar{\omega}_N$ implies that $|\bar{\delta}^*_i|_{\bar{\omega}_N} \le r_i \le \varphi^{-1}(\epsilon)$ for any $i=1,2,\cdots, d$, from property (ii) of Lemma \ref{lem:LgJ}, we get (\ref{eqn:PNomegaJstargPhi}). $\Box$

From Theorem \ref{thm:gNJstar}, we can derive an explicit ULB as follow.

\begin{corollary}
	Suppose the conditions in Theorem \ref{thm:gNJstar} hold, a ULB satisfying (\ref{eqn:pDepsilon}) can be given as follow:
	\begin{align}\label{eqn:tildealphaLSPL}
		\tilde{\alpha}(\varepsilon) = L_g\varphi^{-1}(1-\sqrt[N]{\frac{1-\varepsilon}{d}}), \forall \varepsilon \in (0,1). 
	\end{align}
\end{corollary}
Proof: Letting $\varepsilon = 1-d(1-\epsilon)^N$, this result is a direct consequence of Theorem \ref{thm:gNJstar}. $\Box$

According to Remark \ref{rem:phi}, for any $\epsilon\in (0,1/d]$, we can even replace the bound in (\ref{eqn:BPhi}) with the one in (\ref{eqn:PhiN}) in Theorem \ref{thm:gNJstar}. Thus, a tighter ULB can be derived. We omit the details due to page limitation.



\subsection{Lipschitz continuity}
In the rest of this section, we discuss sufficient conditions for the continuity condition \eqref{eqn:gdLg}.  First, we need Slater's condition of the RCP as given below.
\begin{assumption}\label{ass:slater}
(Slater's condition) There exists $x_0\in X$ such that $\sup_{\delta\in \Delta} f(x_0,\delta) <0$.
\end{assumption}
Under this assumption, following \cite{ART:ESL14}, let us define the following constant:
\begin{align}\label{eqn:LSP}
L_{SP} \coloneqq \frac{\min_{x\in X}c^\top x-c^\top x_0}{\sup_{\delta \in \Delta}f(x_0,\delta)}
\end{align}
We then assume Lipschitz continuity on $f(x,\delta)$ in $\delta$ inside the closure of $\Delta$.

\begin{assumption}\label{ass:Lipschitz}
For any $x\in X$, $f(x,\delta)$ is Lipschitz continuous in $\delta$ with constant $L_{\delta}$ in $\textrm{cl}(\Delta)$, i.e.,
\begin{align}\label{eqn:fLipschitz}
\forall \delta_1,\delta_2 \in \textrm{cl}(\Delta), \|f(x,\delta_1)-f(x,\delta_2)\| \le L_{\delta} \|\delta_1-\delta_2\|, 
\end{align}
where $\textrm{cl}(\Delta)$ denote the closure of $\Delta$.
\end{assumption}

Based on these two additional assumptions, we claim that the continuity condition \eqref{eqn:gdLg} is guaranteed, as shown in the next proposition. 

\begin{proposition}\label{prop:lipschitz}
Given any $N\in \mathbb{Z}^+$, let $g_N(\omega)$ be defined as in (\ref{eqn:gddelta}) for any $\omega\in \textrm{cl}(\Delta)^N$. Suppose Assumptions \ref{ass:convex} -- \ref{ass:Lipschitz} hold, then, for any $\omega=(\delta_1,\delta_2,\cdots,\delta_N) \in \textrm{cl}(\Delta)^N$ and $\omega'=(\delta'_1,\delta'_2,\cdots,\delta'_N) \in \textrm{cl}(\Delta)^N$,  it holds that
\begin{align}\label{eqn:gNLg}
	\|g_N(\omega)-g_N(\omega')\| \le L_{SP}L_{\delta} \max_{i\in \{1,2,\cdots, N\}} \|\delta_i-\delta'_i\|.
\end{align}
where $L_{SP}$ and $L_{\delta}$ are given in (\ref{eqn:LSP}) and (\ref{eqn:fLipschitz}) respectively. 
\end{proposition}
Proof: For any two points $\omega = (\delta_1,\delta_2,\cdots, \delta_N),\omega' = (\delta'_1,\delta'_2,\cdots, \delta'_N)\in \textrm{cl}(\Delta)^N$, we want to show \eqref{eqn:gNLg}. First, we defined the perturbed problem for $\omega$ with the perturbation $L_{\delta} \max_j \|\delta_j-\delta'_j\|$:
\begin{subequations}
\begin{align}
&\tilde{g}_N(\omega) \coloneqq \min_{x\in X} c^\top x\\
\textrm{s.t.}  & \quad f(x,\delta_i) \le L_{\delta} \max_j \|\delta_j-\delta'_j\|, \forall i=1,2,\cdots, N
\end{align}
\end{subequations}
Similarly, the perturbed problem for $\omega'$ is also defined: 
\begin{subequations}
\begin{align}
&\tilde{g}_N(\omega') \coloneqq \min_{x\in X} c^\top x\\
\textrm{s.t.} & \quad f(x,\delta'_i) \le L_{\delta} \max_j \|\delta_j-\delta'_j\|,\forall i=1,2,\cdots, N
\end{align}
\end{subequations}
From strong duality and perturbation analysis, see, e.g., Section 5.6 in \cite{BOO:BV04} and Lemma 3.4 in \cite{ART:ESL14}, the optimal cost function of the perturbed problem is Lipschitz continuous with respect to the perturbation with constant $L_{SP}$. Thus, $0 \le g_N(\omega)-\tilde{g}_N(\omega)\le L_{SP}L_{\delta} \max_j \|\delta_j-\delta'_j\|$ and $0\le g_N(\omega')-\tilde{g}_N(\omega') \le L_{SP}L_{\delta} \max_j \|\delta_j-\delta'_j\|$. From the Lipschitz continuity of $f(x,\delta)$ in $\delta$, $|f(x,\delta_i) -f(x,\delta'_i)| \le L_{\delta}\|\delta_i-\delta'_i\| \le L_{\delta}\max_j \|\delta_j-\delta'_j\|$ for all $i=1,2,\cdots, N$, which implies that  $\{x: f(x,\delta_i) \le 0, i=1,2,\cdots, N\} \subseteq \{x: f(x,\delta'_i) \le L_{\delta} \max_j \|\delta_j-\delta'_j\|, i=1,2,\cdots, N\}$ and $\{x: f(x,\delta'_i) \le 0, i=1,2,\cdots, N\} \subseteq \{x: f(x,\delta_i) \le L_{\delta} \max_j \|\delta_j-\delta'_j\|, i=1,2,\cdots, N\}$. Hence, $\tilde{g}_N(\omega) \le g_N(\omega')$ and  $\tilde{g}_N(\omega') \le g_N(\omega)$. Finally, we conclude that $|g_N(\omega')-g_N(\omega)|\le L_{SP}L_{\delta} \max_j \|\delta_j-\delta'_j\|$. $\Box$

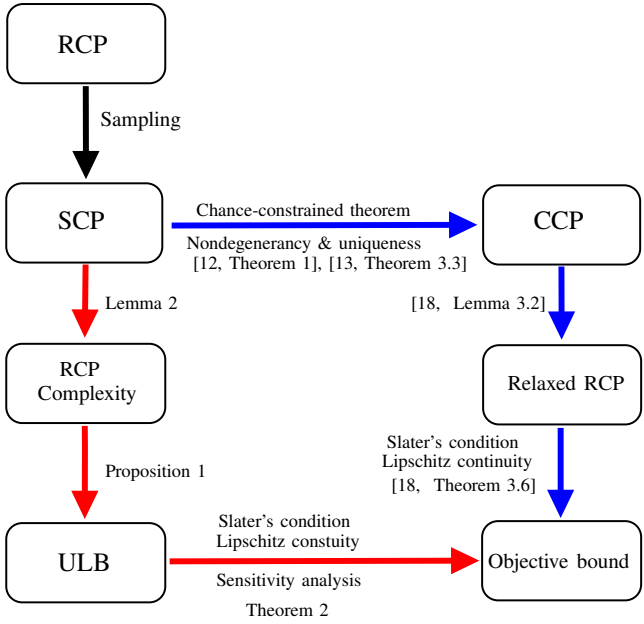
\begin{figure}
\tikzset{every picture/.style={line width=0.75pt}} 

\begin{tikzpicture}[x=0.75pt,y=0.75pt,yscale=-1,xscale=1]
	
	\draw   (111,77) .. controls (111,72.58) and (114.58,69) .. (119,69) -- (182.7,69) .. controls (187.12,69) and (190.7,72.58) .. (190.7,77) -- (190.7,101) .. controls (190.7,105.42) and (187.12,109) .. (182.7,109) -- (119,109) .. controls (114.58,109) and (111,105.42) .. (111,101) -- cycle ;
	\draw   (110.5,167.5) .. controls (110.5,163.08) and (114.08,159.5) .. (118.5,159.5) -- (181.2,159.5) .. controls (185.62,159.5) and (189.2,163.08) .. (189.2,167.5) -- (189.2,191.5) .. controls (189.2,195.92) and (185.62,199.5) .. (181.2,199.5) -- (118.5,199.5) .. controls (114.08,199.5) and (110.5,195.92) .. (110.5,191.5) -- cycle ;
	\draw [color={rgb, 255:red, 0; green, 0; blue, 255 }  ,draw opacity=1 ][line width=2.25]    (193.9,180.8) -- (343.1,180.46) ;
	\draw [shift={(348.1,180.45)}, rotate = 179.87] [fill={rgb, 255:red, 0; green, 0; blue, 255 }  ,fill opacity=1 ][line width=0.08]  [draw opacity=0] (14.29,-6.86) -- (0,0) -- (14.29,6.86) -- cycle    ;
	\draw   (351,168.5) .. controls (351,164.08) and (354.58,160.5) .. (359,160.5) -- (422.2,160.5) .. controls (426.62,160.5) and (430.2,164.08) .. (430.2,168.5) -- (430.2,192.5) .. controls (430.2,196.92) and (426.62,200.5) .. (422.2,200.5) -- (359,200.5) .. controls (354.58,200.5) and (351,196.92) .. (351,192.5) -- cycle ;
	\draw [line width=2.25]    (150.6,110.45) -- (150.63,122.45) -- (150.69,152.45) ;
	\draw [shift={(150.7,157.45)}, rotate = 269.88] [fill={rgb, 255:red, 0; green, 0; blue, 0 }  ][line width=0.08]  [draw opacity=0] (14.29,-6.86) -- (0,0) -- (14.29,6.86) -- cycle    ;
	\draw   (350.5,338) .. controls (350.5,333.58) and (354.08,330) .. (358.5,330) -- (422.7,330) .. controls (427.12,330) and (430.7,333.58) .. (430.7,338) -- (430.7,362) .. controls (430.7,366.42) and (427.12,370) .. (422.7,370) -- (358.5,370) .. controls (354.08,370) and (350.5,366.42) .. (350.5,362) -- cycle ;
	\draw [color={rgb, 255:red, 0; green, 0; blue, 255 }  ,draw opacity=1 ][line width=2.25]    (390.5,201.5) -- (390.59,233.45) ;
	\draw [shift={(390.6,238.45)}, rotate = 269.84] [fill={rgb, 255:red, 0; green, 0; blue, 255 }  ,fill opacity=1 ][line width=0.08]  [draw opacity=0] (14.29,-6.86) -- (0,0) -- (14.29,6.86) -- cycle    ;
	\draw [color={rgb, 255:red, 255; green, 0; blue, 0 }  ,draw opacity=1 ][line width=2.25]    (150.5,200.5) -- (150.36,214.43) -- (150.16,232.45) ;
	\draw [shift={(150.1,237.45)}, rotate = 270.65] [fill={rgb, 255:red, 255; green, 0; blue, 0 }  ,fill opacity=1 ][line width=0.08]  [draw opacity=0] (14.29,-6.86) -- (0,0) -- (14.29,6.86) -- cycle    ;
	\draw   (110.5,339) .. controls (110.5,334.58) and (114.08,331) .. (118.5,331) -- (182.7,331) .. controls (187.12,331) and (190.7,334.58) .. (190.7,339) -- (190.7,363) .. controls (190.7,367.42) and (187.12,371) .. (182.7,371) -- (118.5,371) .. controls (114.08,371) and (110.5,367.42) .. (110.5,363) -- cycle ;
	\draw [color={rgb, 255:red, 255; green, 0; blue, 0 }  ,draw opacity=1 ][line width=2.25]    (192.9,350.3) -- (344.1,349.96) ;
	\draw [shift={(349.1,349.95)}, rotate = 179.87] [fill={rgb, 255:red, 255; green, 0; blue, 0 }  ,fill opacity=1 ][line width=0.08]  [draw opacity=0] (14.29,-6.86) -- (0,0) -- (14.29,6.86) -- cycle    ;
	\draw   (351,249.2) .. controls (351,244.78) and (354.58,241.2) .. (359,241.2) -- (422.7,241.2) .. controls (427.12,241.2) and (430.7,244.78) .. (430.7,249.2) -- (430.7,273.2) .. controls (430.7,277.62) and (427.12,281.2) .. (422.7,281.2) -- (359,281.2) .. controls (354.58,281.2) and (351,277.62) .. (351,273.2) -- cycle ;
	\draw [color={rgb, 255:red, 0; green, 0; blue, 255 }  ,draw opacity=1 ][line width=2.25]    (390.1,282.95) -- (390.1,323.45) ;
	\draw [shift={(390.1,328.45)}, rotate = 270] [fill={rgb, 255:red, 0; green, 0; blue, 255 }  ,fill opacity=1 ][line width=0.08]  [draw opacity=0] (14.29,-6.86) -- (0,0) -- (14.29,6.86) -- cycle    ;
	\draw   (110.5,248.5) .. controls (110.5,244.08) and (114.08,240.5) .. (118.5,240.5) -- (182.7,240.5) .. controls (187.12,240.5) and (190.7,244.08) .. (190.7,248.5) -- (190.7,272.5) .. controls (190.7,276.92) and (187.12,280.5) .. (182.7,280.5) -- (118.5,280.5) .. controls (114.08,280.5) and (110.5,276.92) .. (110.5,272.5) -- cycle ;
	\draw [color={rgb, 255:red, 255; green, 0; blue, 0 }  ,draw opacity=1 ][line width=2.25]    (150.1,281.95) -- (150.18,292.44) -- (150.11,323.95) ;
	\draw [shift={(150.1,328.95)}, rotate = 270.13] [fill={rgb, 255:red, 255; green, 0; blue, 0 }  ,fill opacity=1 ][line width=0.08]  [draw opacity=0] (14.29,-6.86) -- (0,0) -- (14.29,6.86) -- cycle    ;

	\draw (135,83) node [anchor=north west][inner sep=0.75pt]   [align=left] {RCP};
	\draw (135,172) node [anchor=north west][inner sep=0.75pt]   [align=left] {SCP};
	\draw (205,168) node [anchor=north west][inner sep=0.75pt]   [align=left] {{\scriptsize Chance-constrained theorem}};
	\draw (200,185) node [anchor=north west][inner sep=0.75pt]   [align=left] {{\scriptsize Nondegenerancy \& uniqueness}};
	\draw (375,172) node [anchor=north west][inner sep=0.75pt]   [align=left] {CCP};
	\draw (352,345) node [anchor=north west][inner sep=0.75pt]   [align=left] {{\footnotesize Objective bound}};
	\draw (156.5,121.7) node [anchor=north west][inner sep=0.75pt]   [align=left] {{\footnotesize Sampling}};
	\draw (135,345) node [anchor=north west][inner sep=0.75pt]   [align=left] {ULB};
	\draw (159,215) node [anchor=north west][inner sep=0.75pt]   [align=left] {{\scriptsize Lemma \ref{lem:JG}}};
	\draw (216,325) node [anchor=north west][inner sep=0.75pt]   [align=left] {{\scriptsize Slater's condition}};
	\draw (215,335) node [anchor=north west][inner sep=0.75pt]   [align=left] {{\scriptsize Lipschitz constuity}};
	\draw (362,255) node [anchor=north west][inner sep=0.75pt]   [align=left] {{\footnotesize Relaxed RCP}};
	\draw (301,285) node [anchor=north west][inner sep=0.75pt]   [align=left] {{\scriptsize Slater's condition}};
	\draw (298.5,295) node [anchor=north west][inner sep=0.75pt]   [align=left] {{\scriptsize Lipschitz continuity}};
	\draw (136,248) node [anchor=north west][inner sep=0.75pt]   [align=left] {{\footnotesize RCP }};
	\draw (125,260) node [anchor=north west][inner sep=0.75pt]   [align=left] {{\footnotesize Complexity}};
	\draw (215,355) node [anchor=north west][inner sep=0.75pt]   [align=left] {{\scriptsize Sensitivity analysis}};
	\draw (310,215) node [anchor=north west][inner sep=0.75pt]   [align=left] {{\scriptsize \cite[ Lemma 3.2]{ART:ESL14}}};
	\draw (200,195) node [anchor=north west][inner sep=0.75pt]   [align=left] {{\scriptsize \cite[Theorem 1]{ART:CG08},\cite[Theorem 3.3]{ART:C10}}};
	\draw (300,307) node [anchor=north west][inner sep=0.75pt]   [align=left] {{\scriptsize \cite[ Theorem 3.6]{ART:ESL14}}};
	\draw (159,300) node [anchor=north west][inner sep=0.75pt]   [align=left] {{\scriptsize Proposition \ref{prop:boundd}}};
	\draw (230,370) node [anchor=north west][inner sep=0.75pt]   [align=left] {{\scriptsize Theorem \ref{thm:gNJstar}}};
\end{tikzpicture}
\caption{Flowchart of the proposed method and its comparison with the two-step approach: The red arrows represent our method and the blue arrows represent the two-step approach.}\label{fig:sum}
\end{figure}

The continuity property above then allows to derive explicit ULBs as mentioned above. The overall procedure is summarized in Figure \ref{fig:sum}. Putting Theorem \ref{thm:gNJstar} and Proposition \ref{prop:lipschitz} together, we immediately arrive at the following corollary.

\begin{corollary}
	Consider the same conditions in Theorem \ref{thm:gNJstar} and suppose Assumptions \ref{ass:slater} \& \ref{ass:Lipschitz} hold.  For any $\beta\in (0,1)$, it holds that
	\begin{align}
		\mathbb{P}^N\{\omega_N \in \Delta^{N}: &J^\star -g_N(\omega_N) \\
		\le &L_{SP}L_{\delta} \varphi^{-1}(\Phi_a^{-1}(\beta;d,N))\} \ge 1-\beta. \nonumber
	\end{align}
where $\Phi_a^{-1}(\cdot;d,N)$ denote the inverse function of $\Phi_a(\cdot;d,N)$ with fixed $d$ and $N$.
\end{corollary}
Proof: This is a direct consequence of Theorem \ref{thm:gNJstar} and Proposition \ref{prop:lipschitz}. $\Box$

\begin{remark}
When $\beta \ge \Phi(1/d;d,N)$, as mentioned in Remark \ref{rem:phi}, a better bound can be obtained
\begin{align}
	\mathbb{P}^N\{\omega_N \in \Delta^{N}: &J^\star -g_N(\omega_N) \\
	\le &L_{SP}L_{\delta} \varphi^{-1}(\Phi^{-1}(\beta;d,N))\} \ge 1-\beta. \nonumber
\end{align}
As a comparison, we also show the bound from Theorem 3.6 in \cite{ART:ESL14} below
\begin{align}
	\mathbb{P}^N\{\omega_N \in \Delta^{N}: &J^\star -g_N(\omega_N) \\
	\le &L_{SP}L_{\delta} \varphi^{-1}(\Phi^{-1}_c(\beta;d,N))\} \ge 1-\beta. \nonumber
\end{align}
The difference of these bounds lies in the probability bounds $\Phi_a(\beta;d,N)$, $\Phi(\beta;d,N)$ and $\Phi_c(\beta;d,N)$, which are depicted in Figure \ref{fig:phi}.
\end{remark}


\begin{remark}
	We also want to mention that  the proposed method can be extended to non-convex programs in which the complexity (see Definition \ref{def:helly}) or its upper bound is known a priori with proper continuity and regularity conditions. We leave the details for future work.
\end{remark}

\section{An illustrative example}\label{sec:example}
\begin{figure}
\centering

\tikzset{every picture/.style={line width=0.5pt}} 

\begin{tikzpicture}[x=0.5pt,y=0.5pt,yscale=-1,xscale=1]
	
	\draw  (153,191) -- (505,191)(318.55,42) -- (318.55,313) (498,186) -- (505,191) -- (498,196) (313.55,49) -- (318.55,42) -- (323.55,49)  ;
	\draw   (244.25,191) .. controls (244.25,149.72) and (277.72,116.25) .. (319,116.25) .. controls (360.28,116.25) and (393.75,149.72) .. (393.75,191) .. controls (393.75,232.28) and (360.28,265.75) .. (319,265.75) .. controls (277.72,265.75) and (244.25,232.28) .. (244.25,191) -- cycle ;
	\draw  [dash pattern={on 4.5pt off 4.5pt}]  (308,86) -- (378.52,155.49) -- (444,220) ;
	\draw  [dash pattern={on 4.5pt off 4.5pt}]  (205,171) -- (275.52,240.49) -- (341,305) ;
	\draw    (372,137) -- (266,244) ;
	\draw [shift={(266,244)}, rotate = 134.73] [color={rgb, 255:red, 0; green, 0; blue, 0 }  ][fill={rgb, 255:red, 0; green, 0; blue, 0 }  ][line width=0.75]      (0, 0) circle [x radius= 3.35, y radius= 3.35]   ;
	\draw [shift={(319,190.5)}, rotate = 134.73] [color={rgb, 255:red, 0; green, 0; blue, 0 }  ][fill={rgb, 255:red, 0; green, 0; blue, 0 }  ][line width=0.75]      (0, 0) circle [x radius= 3.35, y radius= 3.35]   ;
	\draw [shift={(372,137)}, rotate = 134.73] [color={rgb, 255:red, 0; green, 0; blue, 0 }  ][fill={rgb, 255:red, 0; green, 0; blue, 0 }  ][line width=0.75]      (0, 0) circle [x radius= 3.35, y radius= 3.35]   ;
	\draw  [draw opacity=0][line width=1.5]  (341.78,119.08) .. controls (364,125.72) and (381.77,142.57) .. (389.52,164.14) -- (320.65,188.5) -- cycle ; \draw  [color={rgb, 255:red, 255; green, 0; blue, 0 }  ,draw opacity=1 ][line width=1.5]  (341.78,119.08) .. controls (364,125.72) and (381.77,142.57) .. (389.52,164.14) ;  
	\draw  [draw opacity=0][line width=1.5]  (297.21,261.2) .. controls (272.88,254.97) and (253.46,236.62) .. (245.88,213.03) -- (315.55,191) -- cycle ; \draw  [color={rgb, 255:red, 0; green, 0; blue, 255 }  ,draw opacity=1 ][line width=1.5]  (297.21,261.2) .. controls (272.88,254.97) and (253.46,236.62) .. (245.88,213.03) ;  
	
	\draw (305.5,172.4) node [anchor=north west][inner sep=0.75pt]    {$0$};
	\draw (395.75,194.4) node [anchor=north west][inner sep=0.75pt]    {$1$};
	\draw (303.5,94.4) node [anchor=north west][inner sep=0.75pt]    {$1$};
	\draw (376,128.4) node [anchor=north west][inner sep=0.75pt]  [font=\footnotesize]  {$\left(\sqrt{2} ,\sqrt{2}\right)$};
	\draw (210,262.4) node [anchor=north west][inner sep=0.75pt]  [font=\footnotesize]  {$\left( -\sqrt{2} ,-\sqrt{2}\right)$};
	\draw (473,201.4) node [anchor=north west][inner sep=0.75pt]    {$\delta _{1}$};
	\draw (294,44.4) node [anchor=north west][inner sep=0.75pt]    {$\delta _{2}$};
	\draw (415,227.4) node [anchor=north west][inner sep=0.75pt]  [font=\footnotesize]  {$b^{\top } \delta =\sqrt{2} -\alpha $};
	\draw (354,102.4) node [anchor=north west][inner sep=0.75pt]  [font=\footnotesize]  {$\overline{S}_{\alpha }$};
	\draw (230,220.4) node [anchor=north west][inner sep=0.75pt]  [font=\footnotesize]  {$\underline{S}_{\alpha }$};
	\draw (341,287.4) node [anchor=north west][inner sep=0.75pt]  [font=\footnotesize]  {$b^{\top } \delta =\alpha -\sqrt{2}$};

\end{tikzpicture}
\caption{Illustration of $\overline{S}_\alpha$ and$ \underline{S}_\alpha$.}\label{fig:exS}
\end{figure}
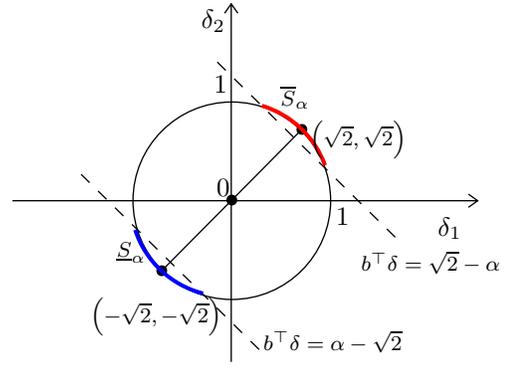

\begin{figure}
	\centering
	\includegraphics[width=0.9\linewidth]{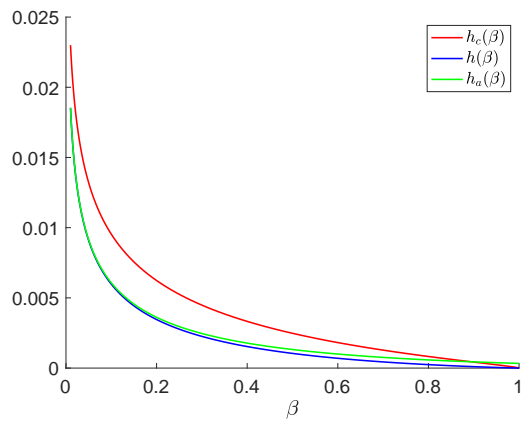}
	\caption{Comparison of different objective bounds.}\label{fig:ex}
\end{figure}

In this section, we present an example to illustrate the conservatism of the two-step approach in \cite{ART:ESL14} induced by the chance-constrained theorems in \cite{ART:CG08,ART:C10}. Consider the following robust optimization problem
\begin{align*}
	&\min_{x_1,x_2} x_2\\
	\textrm{s.t.} \quad &|x_1+ b^\top \delta| \le x_2, \forall \delta \in \Delta = \{\delta\in \mathbb{R}^2: \|\delta\| =1\},   \\
	& |x_1| \le 1, |x_2| \le \sqrt{2}.
\end{align*}
where $b = [1;1]$. As $b^\top \delta \in [-\sqrt{2},\sqrt{2}]$ for all $\|\delta\|=1$, the optimal $J^\star = \sqrt{2}$ with the optimizer being $x_1^\star =0$ and $x_2^\star = \sqrt{2}$. Given $N$ i.i.d. samples $\omega_N \coloneqq \{\delta_1,\delta_2,\cdots, \delta_N\}$ for some $N\in \mathbb{Z}^+$, drawn from the uniform distribution $\mathbb{P}$ on  $\Delta$, let
\begin{align*}
	\underline{\eta}(\omega_N) \coloneqq \min_{i=1,2,\cdots, N} b^\top  \delta_i, \quad \overline{\eta}(\omega_N) \coloneqq \max_{i=1,2,\cdots,N} b^\top \delta_i.
	\end{align*}
With these definitions, the optimizer of the sampled problem SCP is 
\begin{align*}
	x^*_1(\omega_N) &= \frac{\overline{\eta}(\omega_N)+\underline{\eta}(\omega_N)}{2}, x^*_2(\omega_N) = \frac{\overline{\eta}(\omega_N)-\underline{\eta}(\omega_N)}{2}.
\end{align*}
The tail probability $p_N(\alpha)$ defined as in (\ref{eqn:palpha}) becomes 
\begin{align*}
	p_N(\alpha) = \mathbb{P}^N\{\omega_N \in \Delta^N: \sqrt{2}-\alpha \le \frac{\overline{\eta}(\omega_N)-\underline{\eta}(\omega_N)}{2}\} 
\end{align*}
for $\alpha \in [0,\sqrt{2}]$. Given any $\alpha \in [0,\sqrt{2}]$, a relaxation of $p_N(\alpha)$ is given below
\begin{align*}
	\tilde{p}_N(\alpha) = \mathbb{P}^N\{\omega_N \in \Delta^N: &\overline{\eta}(\omega_N) \ge \sqrt{2}-\alpha, \\
	&\underline{\eta}(\omega_N) \le \alpha-\sqrt{2}\}.
\end{align*}
We also define the following sets:
$
\overline{S}_\alpha \coloneqq \{\delta\in \Delta: b^\top \delta \ge \sqrt{2}-\alpha\}
$
 and 
$
\underline{S}_\alpha \coloneqq \{\delta\in \Delta: b^\top \delta \le \alpha-\sqrt{2}\}
$, as illustrated in Figure \ref{fig:exS}. With simple manipulations, 
$
\mathbb{P}\{ \overline{S}_\alpha  \} = \mathbb{P}\{ \underline{S}_\alpha  \}= \frac{\cos^{-1}(\frac{\sqrt{2}-\alpha}{\sqrt{2}})}{\pi}
$.
The relaxation $\tilde{p}_N(\alpha)$ can be considered as the probability that $\bar{\omega}_N \cap \overline{S}_\alpha \not= \emptyset$ and $\bar{\omega}_N \cap \underline{S}_\alpha \not= \emptyset$ with $\bar{\omega}_N\coloneqq \{\delta_1,\delta_2,\cdots, \delta_N\}$. Hence, from Lemma \ref{lem:Phi} and Remark \ref{rem:phi}, 
\begin{align*}
		\tilde{p}_N(\alpha) &= 1-\Phi(\frac{\cos^{-1}(\frac{\sqrt{2}-\alpha}{\sqrt{2}})}{\pi};2,N)\\
		& \le 1-\Phi_a(\frac{\cos^{-1}(\frac{\sqrt{2}-\alpha}{\sqrt{2}})}{\pi};k,N)
\end{align*}
Now, let us consider the optimal ULB of $\tilde{p}_N(\alpha)$:
$
\tilde{\alpha}^*(\varepsilon)\coloneqq \sup_{\alpha\ge 0} \alpha: \tilde{p}_N(\alpha)\le \varepsilon, \forall \varepsilon\in (0,1). 
$
As $\tilde{p}_N(\alpha) \le p_N(\alpha)$, 
\begin{align*}
\alpha^*(\varepsilon) \le \tilde{\alpha}^*(\varepsilon) =  h(1-\varepsilon) \le h_a(1-\varepsilon)
\end{align*}
where $h_a(\beta)\coloneqq \sqrt{2}(1-\cos(\pi \Phi_a^{-1}(\beta;2,N) ))$ and  $h(\beta)\coloneqq  \sqrt{2}(1-\cos(\pi \Phi^{-1}(\beta;2,N) ))$, $\beta \in (0,1)$.  Finally, we conclude that, given the confidence level $\beta\in (0,1)$
\begin{align}
		&\mathbb{P}^{N} \{\omega_N \in \Delta^N: J^\star-x^*_2(\omega_N) \le h_a(\beta)  \} \nonumber \\
		\ge & \mathbb{P}^{N} \{\omega_N \in \Delta^N: J^\star-x^*_2(\omega_N) \le h(\beta)  \} 
	\ge 1- \beta.
\end{align}
We now consider Theorem \ref{thm:chance}. For any $\epsilon\in (0,1)$, 
\begin{align*}
	&\mathbb{P}^{N} \{\omega_N \in \Delta^N: \mathbb{P}\{\delta: J(\bar{\omega}_N \cup \{\delta\}) > J(\bar{\omega}_N)\} \le  \epsilon \} \\
	\ge &1- \Phi_c(\epsilon;d,N)
\end{align*}
In fact, the equality holds because this is a fully supported problem as defined in \cite{ART:C10}. We can see that
\begin{align*}
	&\{\delta: J(\bar{\omega}_N \cup \{\delta\}) > J(\bar{\omega}_N)\} \\
	=& \{\delta: b^\top \delta > \overline{\eta}(\omega_N)\} \cup \{\delta: b^\top \delta < \underline{\eta}(\omega_N)\}.
\end{align*}
Hence, $\mathbb{P}\{\delta: J(\bar{\omega}_N \cup \{\delta\}) > J(\bar{\omega}_N)\} \le  \epsilon$ implies 
\begin{align}\label{eqn:Pepsilon}
	\mathbb{P}\{\delta: b^\top \delta > \overline{\eta}(\omega_N)\}  + \mathbb{P} \{\delta: b^\top \delta < \underline{\eta}(\omega_N)\} \le \epsilon.
\end{align}
Then, we need to derive the objective bound for any case. We find out that the worst case is when  $	\mathbb{P}\{\delta: b^\top \delta > \overline{\eta}(\omega_N)\}  = \epsilon$ and $ \mathbb{P} \{\delta: b^\top \delta < \underline{\eta}(\omega_N)\}=0$ (or vice versa), which means that $\overline{\eta}(\omega_N) = \sqrt{2}\cos(\pi\epsilon) $  and $\underline{\eta}(\omega_N) =-\sqrt{2}$. In this case, $x^*_2(\omega_N)  = \frac{\sqrt{2}(1+\cos(\pi \epsilon))}{2}$. Hence, the objective value bound we can claim is 
\begin{align*}
	&\mathbb{P}^{N} \{\omega_N \in \Delta^N: J^\star-x^*_2(\omega_N) \le \frac{\sqrt{2}(1-\cos(\pi \epsilon))}{2} \} \\
	\ge &1- \Phi_c(\epsilon;d,N).
\end{align*}
For any confidence level $\beta\in (0,1)$, 
\begin{align}
	\mathbb{P}^{N} \{\omega_N \in \Delta^N: J^\star-x^*_2(\omega_N) \le h_c(\beta)  \} \ge 1-\beta,
\end{align}
where $h_c(\beta)\coloneqq \frac{\sqrt{2}(1-\cos(\pi \Phi_c^{-1}(\beta;2,N) ))}{2}$. The comparison of these objective bounds is shown in Figure \ref{fig:ex} with $N = 100$. From this example, we can see that Theorem \ref{thm:chance} induces conservatism into the derivation of the objective value bound as the worst case in \eqref{eqn:Pepsilon} has to be taken into account.

\section{Conclusions}
In this paper, we have presented a novel method for deriving probabilistic guarantees on the objective value for the scenario approach for robust optimization problems via sensitivity analysis. This method is based in particular on the concept of \emph{complexity} of a robust optimization problem, which is the minimal number of constraints that define the actual solution. For robust convex programs, we show that the \emph{complexity} is bounded by the dimension of the decision variable by the use of a max-min reformulation and thereby extend classical results, which are known for scenario programs but not the original robust optimization problem, which possibly contains an infinite number of constraints. With this observation and additional regularity conditions, we are able to build probabilistic objective value bounds for the scenario approach. Our results do not rely on the chance-constrained theorems in the literature as an intermediate step, which enables us to release some strong technical assumptions for ensuring these chance-constrained theorems. Under the same conditions, our bounds also outperform an existing result based on a tight chance-constrained theorem. We illustrate these improvements using a numerical example.

\bibliographystyle{unsrt}
\bibliography{Reference}

\begin{thebibliography}{10}

\bibitem{BOO:P95}
A.~Pr{\'e}kopa.
\newblock {\em Stochastic programming}, volume 324.
\newblock Springer Science \& Business Media, 1995.

\bibitem{ART:RS03}
A.~Ruszczy{\'n}ski and A.~Shapiro.
\newblock Stochastic programming models.
\newblock {\em Handbooks in operations research and management science},
  10:1--64, 2003.

\bibitem{BOO:BL11}
J.~R. Birge and F.~Louveaux.
\newblock {\em Introduction to stochastic programming}.
\newblock Springer Science \& Business Media, 2011.

\bibitem{ART:BN98}
A.~Ben-Tal and A.~Nemirovski.
\newblock Robust convex optimization.
\newblock {\em Mathematics of Operations Research}, 23(4):769--805, 1998.

\bibitem{ART:BN08}
Aharon Ben-Tal and Arkadi Nemirovski.
\newblock Selected topics in robust convex optimization.
\newblock {\em Mathematical Programming}, 112(1):125--158, 2008.

\bibitem{BOO:BEN09}
A.~Ben-Tal, L.~El~Ghaoui, and A.~Nemirovski.
\newblock {\em Robust optimization}.
\newblock Princeton university press, 2009.

\bibitem{ART:BBC11}
Dimitris Bertsimas, David~B Brown, and Constantine Caramanis.
\newblock Theory and applications of robust optimization.
\newblock {\em SIAM review}, 53(3):464--501, 2011.

\bibitem{ART:CC05}
G.~Calafiore and M.~C. Campi.
\newblock Uncertain convex programs: randomized solutions and confidence
  levels.
\newblock {\em Mathematical Programming}, 102(1):25--46, 2005.

\bibitem{ART:CC06}
G.~Calafiore and M.~C. Campi.
\newblock The scenario approach to robust control design.
\newblock {\em IEEE Transactions on Automatic Control}, 51(5):742--753, 2006.

\bibitem{BOO:CG18}
M.~C. Campi and S.~Garatti.
\newblock {\em Introduction to the scenario approach}.
\newblock SIAM, 2018.

\bibitem{ART:CCG21}
M.C. Campi, A.~Car{\`e}, and S.~Garatti.
\newblock The scenario approach: A tool at the service of data-driven decision
  making.
\newblock {\em Annual Reviews in Control}, 52:1--17, 2021.

\bibitem{ART:CG08}
M.~C. Campi and S.~Garatti.
\newblock The exact feasibility of randomized solutions of uncertain convex
  programs.
\newblock {\em SIAM Journal on Optimization}, 19(3):1211--1230, 2008.

\bibitem{ART:C10}
G.~C. Calafiore.
\newblock Random convex programs.
\newblock {\em SIAM Journal on Optimization}, 20(6):3427--3464, 2010.

\bibitem{ART:CG18}
M.~C. Campi and S.~Garatti.
\newblock Wait-and-judge scenario optimization.
\newblock {\em Mathematical Programming}, 167(1):155--189, 2018.

\bibitem{ART:GC19}
S.~Garatti and M.C. Campi.
\newblock Risk and complexity in scenario optimization.
\newblock {\em Mathematical Programming}, pages 1--37, 2019.

\bibitem{ART:CGR18}
M.~C. Campi, S.~Garatti, and F.~A. Ramponi.
\newblock A general scenario theory for nonconvex optimization and decision
  making.
\newblock {\em IEEE Transactions on Automatic Control}, 63(12):4067--4078,
  2018.

\bibitem{ART:GC21}
S.~Garatti and M.~C. Campi.
\newblock The risk of making decisions from data through the lens of the
  scenario approach.
\newblock {\em IFAC-PapersOnLine}, 54(7):607--612, 2021.

\bibitem{ART:ESL14}
P.~M. Esfahani, T.~Sutter, and J.~Lygeros.
\newblock Performance bounds for the scenario approach and an extension to a
  class of non-convex programs.
\newblock {\em IEEE Transactions on Automatic Control}, 60(1):46--58, 2014.

\bibitem{INP:KDSA14}
J.~Kapinski, J.~V. Deshmukh, S.~Sankaranarayanan, and N.~Arechiga.
\newblock Simulation-guided lyapunov analysis for hybrid dynamical systems.
\newblock In {\em Proceedings of the 17th international conference on Hybrid
  systems: computation and control}, pages 133--142. ACM, 2014.

\bibitem{INP:KQHT16}
A.~Kozarev, J.~Quindlen, J.~How, and U.~Topcu.
\newblock Case studies in data-driven verification of dynamical systems.
\newblock In {\em Proceedings of the 19th International Conference on Hybrid
  Systems: Computation and Control}, pages 81--86. ACM, 2016.

\bibitem{ART:BTSS20}
N.~M. Boffi, S.~Tu, N.~Matni, J.-J.~E. Slotine, and V.~Sindhwani.
\newblock Learning stability certificates from data.
\newblock {\em arXiv preprint arXiv:2008.05952}, 2020.

\bibitem{ART:MXM20}
G.~Mamakoukas, O.~Xherija, and T.~Murphey.
\newblock Memory-efficient learning of stable linear dynamical systems for
  prediction and control.
\newblock {\em Advances in Neural Information Processing Systems},
  33:13527--13538, 2020.

\bibitem{ART:KBJT19}
J.~Kenanian, A.~Balkan, R.~M. Jungers, and P.~Tabuada.
\newblock Data driven stability analysis of black-box switched linear systems.
\newblock {\em Automatica}, 109:108533, 2019.

\bibitem{ART:WJ20}
Z.~Wang and R.~M. Jungers.
\newblock Scenario-based set invariance verification for black-box nonlinear
  systems.
\newblock {\em IEEE Control Systems Letters}, 5(1):193--198, 2021.

\bibitem{INP:BJW21}
G.~O. Berger, R.~M. Jungers, and Z.~Wang.
\newblock Chance-constrained quasi-convex optimization with application to
  data-driven switched systems control.
\newblock In {\em Learning for Dynamics and Control}, pages 571--583. PMLR,
  2021.

\bibitem{INP:RWJ21}
A.~Rubbens, Z.~Wang, and R.~M. Jungers.
\newblock Data-driven stability analysis of switched linear systems with sum of
  squares guarantees.
\newblock In {\em The 7th IFAC Conference on Analysis and Design of Hybrid
  Systems}. IFAC, 2021.

\bibitem{ART:WJ21b}
Z.~Wang and R.~M. Jungers.
\newblock A data-driven method for computing polyhedral invariant sets of
  black-box switched linear systems.
\newblock {\em IEEE Control Systems Letters}, 5(5):1843 -- 1848, 2021.

\bibitem{BOO:BGKKT83}
B.~Bank, J.~Guddat, D.~Klatte, B.~Kummer, and K.~Tammer.
\newblock {\em Non-linear parametric optimization}.
\newblock Birkh{\"a}user Verlag, Basel, 1983.

\bibitem{BOO:RW98}
R.~T. Rockafellar and R.~J-B Wets.
\newblock {\em Variational analysis}, volume 317.
\newblock Springer Science \& Business Media, 1998.

\bibitem{ART:S18}
G.~Still.
\newblock Lectures on parametric optimization: An introduction.
\newblock {\em Optimization Online}, 2018.

\bibitem{BOO:R70}
R.~T. Rockafellar.
\newblock {\em Convex analysis}.
\newblock Princeton university press, 1970.

\bibitem{ART:KT12}
T.~Kanamori and A.~Takeda.
\newblock Worst-case violation of sampled convex programs for optimization with
  uncertainty.
\newblock {\em Journal of Optimization Theory and Applications},
  152(1):171--197, 2012.

\bibitem{ART:CR04}
Antonio Cuevas and Alberto Rodr{\'\i}guez-Casal.
\newblock On boundary estimation.
\newblock {\em Advances in Applied Probability}, 36(2):340--354, 2004.

\bibitem{BOO:BV04}
S.~Boyd and L.~Vandenberghe.
\newblock {\em Convex Optimization}.
\newblock Cambridge University Press, 2004.

\end{thebibliography}

\end{document}